\documentclass[11pt]{amsart}

\usepackage{amsmath,amssymb,epsf, epic, eepic}
\usepackage{amsfonts}
\usepackage{times}

\input xy
\xyoption{all}

\usepackage{slashbox}

\newtheorem{thm}{Theorem}[section]

\newtheorem{prop}[thm]{Proposition}
\newtheorem{exe}[thm]{Example}

\theoremstyle{definition}
\newtheorem{defn}[thm]{Definition}

\newcommand{\bZ}{\mathbb{Z}}
\newcommand{\bQ}{\mathbb{Q}}

\newcommand{\bF}{\mathbb{F}}
\newcommand{\bK}{\mathbb{K}}
\newcommand{\bKprime}{\mathbb{K}^\prime}

\newcommand{\tildeh}{\tilde h}
\newcommand{\tildet}{\tilde t}

\DeclareMathOperator{\id}{Id}

\DeclareMathOperator{\Tor}{Tor}

\newcommand{\ra}{\rightarrow}

\newcommand{\ot}{\otimes}
\newcommand{\ol}{\overline}

\newcommand{\FF}{\bF_2}

\newcommand{\khk}[2]{K\!H^{{#1}}({#2};\bK)}
\newcommand{\uhtk}[2]{U_{h,t}^{{#1}}({#2};\bK)}

\newcommand{\uhtz}[2]{U_{h,t}^{{#1}}({#2};\bZ)}

\newcommand{\uhtq}[2]{U_{h,t}^{{#1}}({#2};\bQ)}

\newcommand{\uhtp}[2]{U_{h,t}^{{#1}}({#2};\bF_p)}
\newcommand{\uhtkt}[2]{U_{\tildeh,\tildet}^{{#1}}({#2};\bK)}

\newcommand{\rasleea}[2]{s({#1},{#2})}
\newcommand{\rasleeknot}[1]{s(K,{#1})}

\newcommand{\rasleemin}[1]{s_{\scriptstyle\text{min}}(K,{#1})}

\newcommand{\rasleemax}[1]{s_{\scriptstyle\text{max}}(K,{#1})}

\title{A remark on Rasmussen's invariant of knots}
\author{Marco Mackaay}
\address{Departamento de Matem\'{a}tica\\ Universidade do Algarve\\ 
Campus de Gambelas\\ 8005-139 Faro\\ Portugal}
\email{mmackaay@ualg.pt}
\author{Paul Turner}
\address{School of Mathematical and Computer Sciences \\Heriot-Watt
  University\\ Edinburgh EH14 4AS\\Scotland}
\email{paul@ma.hw.ac.uk}
\author{Pedro Vaz}
\address{Departamento de Matem\'{a}tica\\ Universidade do Algarve\\ 
Campus de Gambelas\\ 8005-139 Faro\\ Portugal}
\email{pfortevaz@ualg.pt}

\begin{document}


\begin{abstract}
We show that Rasmussen's invariant of knots, which is derived from
Lee's variant of Khovanov homology, is equal to an analogous invariant
derived from certain other filtered link homologies.
\end{abstract}

\maketitle


\section*{Erratum: added 29 June 2012}
We are grateful to  Robert Lipshitz and Sucharit Sarkar for pointing out 
two errors in the proof of Proposition 3.2. in this paper. Firstly, the element $v$, defined in the penultimate displayed 
equation, need not be a cycle over $\mathbb{Z}$ and secondly, the claim that 
$s(\lambda\alpha,\mathbb{Z})=s(\alpha,\mathbb{Z})$, 
just before the last displayed equation, is false in general. 

In consequence the proof of Proposition 3.2 no 
longer holds and the proof of Theorem 4.2,
which relies on it, is no longer valid. The claim in the statement of Theorem 4.2, namely that the Rasmussen invariants defined over $\mathbb Q$ and $\mathbb F_p$ for any prime $p$ are all equal, must, for now, again be considered an open question.

\section{Introduction} 
In \cite{khovanov1} Khovanov introduced a completely new way to define
link invariants. He associated a bigraded cochain complex to a given
link and if two links are ambient isotopic, then the associated
complexes are homotopy equivalent. Thus by taking homology a link
invariant is defined. One of the first variations on Khovanov's
construction was the theory defined by Lee~\cite{lee}. Her link
homology, originally defined over $\bQ$, 
is not bigraded but singly graded with a filtration in place
of what was the internal degree in Khovanov's theory. If one forgets
about the filtration, then Lee's link homology is completely
determined by the linking matrix of the link, which makes it a rather
poor invariant compared to Khovanov's theory. However, by using the
filtration Rasmussen \cite{rasmussen} has defined an integer invariant
of knots $s(K)$ which has many wonderful properties.  For example he
showed that the $s$-invariant yields a lower bound of the smooth slice
genus which led to a new and completely combinatorial proof of the
Milnor conjecture concerning the slice genus of torus knots.  Another
consequence is that if the $s$-invariant of a knot is greater than
zero, then the knot is not smoothly slice which is particularly
interesting if the knot is already known to be topologically
slice. The $s$-invariant is also conjecturally related to the
$\tau$-invariant in Heegaard-Floer knot homology. Much of this is
explained in the survey paper \cite{rasmussensurvey}.

In \cite{barnatan2} Bar-Natan introduced a new link homology theory
defined over $\FF[H]$ where $H$ has internal degree $-2$. Setting
$H=1$ defines a singly graded theory which can be explicitly computed
(see \cite{turner}) and like Lee's theory depends only on the linking
matrix. This theory is again filtered and one can use Rasmussen's
definitions to produce an analogous $s$-invariant using this
theory. The question that motivated the current note was: is
Rasmussen's original $s$-invariant defined using Lee theory the same
as the $s$-invariant defined using Bar-Natan theory? In fact working
over $\bQ$ or $\bF_p$, $p$ a prime, one can define a family of link
homology theories depending on two elements $h$ and $t$, encompassing
Lee's theory and Bar-Natan's theory. Many of these theories give for a
knot a two dimensional vector space in degree zero and for such a
theory one can define a Rasmussen-type invariant.

In the Section \ref{sec:flht} we define the family of link homology
theories of interest to us. We choose the ground field $\bK$ to be one of $\bQ$ or
$\bF_p$, $p$ a prime and the family depends on two parameters $h,t\in\bK$.
We present a couple of computational results and discuss integral
theories. In Section \ref{sec:rsg} we recall Rasmussen's $s$-grading and show
that this is preserved by twist equivalence of theories and by the
universal coefficient theorem. In Section \ref{sec:ri} we define Rasmussen'
$s$-invariant ${\rasleeknot\bK}_{h,t}$ for any theory arising from a
triple $(\bK,h,t)$ such that $h^2+4t$ is a non-zero square in
$\bK$. Letting $\widetilde \bK$ be $\bQ$ or $\bF_p$ ($\bK$ and
$\widetilde\bK$ possibly different) our main result is as follows.

\vspace*{6pt}

{\bf Theorem \ref{thm:main}} 
Let $K$ be a knot. Let $h,t,\tilde h,\tilde t\in \bZ$ be such that
$h^2+4t=\gamma^2\ne 0$ and $\tildeh^2+\tildet=\tilde \gamma^2\ne 0$
with $\gamma\ne 0\in \bK$ and $\tilde\gamma\ne 0\in
\widetilde\bK$. Then
\[
{\rasleeknot\bK}_{h,t} = {\rasleeknot{\widetilde\bK}}_{\tilde h,\tilde t}.
\]


\section{A family of link homology theories}\label{sec:flht}

Let $p$ be a prime and let $\bK$ be $\bQ$ or $\bF_p$. Recall that a {\em Frobenius
system} over $\bK$ is a quadruple $(A,\iota,\Delta,\epsilon)$, where $A$
is a commutative ring with unit $1$, $\iota\colon \bK\to A$ a unital
injective ring homomorphism, $\Delta\colon A\to A\ot A$ a
cocommutative coassociative $A$-bimodule map and $\epsilon\colon A\to
\bK$ a $\bK$-linear map satisfying the additional condition $$(\epsilon\ot
\id)\Delta=\id.$$

Khovanov has explained in ~\cite{khovanov2} how a rank two Frobenius
system gives rise to a link homology theory and moreover that
isomorphic Frobenius systems give rise to isomorphic link homology theories.

\begin{exe} \label{ex:frobsys}
Let $h,t\in \bK$ and define
$$A_{h,t}=\bK[x]/(x^2-hx - t)$$
with coproduct and counit defined by 
$$
\begin{array}{ll}
\Delta(1)=1\ot x+x\ot 1 -h 1\otimes 1,&\Delta(x)=x\ot x+t1\ot 1\\
\epsilon(1)=0,&\epsilon(x)=1.
\end{array}
$$
\end{exe}
This is a rank two Frobenius system which in general is not bi-graded but has a filtration obtained by taking filtration degrees
$\deg(x)=-1$ and $\deg(1)=1$. This filtration induces a filtration on the associated link homology theory. Note that throughout we prefer to use
the grading conventions in \cite{khovanov1} rather than those in
\cite{khovanov2}. These theories are obtained from Khovanov's theory
$A_5$ in \cite{khovanov2} by specialisation of the variable $h$ and
$t$ to elements of $\bK$.  When $h=t=0$ the resulting theory is
Khovanov's original link homology with coefficients in $\bK$ which we
denote $\khk *  -$. In this case the theory is genuinely
bi-graded. When $\bK=\bQ$, $h=0$ and $t=1$ one gets Lee's theory
\cite{lee} and when $\bK=\FF$, $h=1$ and $t=0$ one gets Bar-Natan's
theory \cite{barnatan2}. We will denote the theory defined from
$h,t\in \bK$ by $\uhtk * L$ for a link $L$.

There is one further idea from \cite{khovanov2} that is important for us.
Let $A$ be a Frobenius system and let $\theta\in A$ be an invertible element. 
Then we can {\em twist} $A$ {\em by} $\theta$ to obtain a new
Frobenius system with the same product and unit map but a new
coproduct and counit map defined by $\Delta'(a)=\Delta(\theta^{-1}a)$
and $\epsilon'(a)=\epsilon(\theta a)$.  We call two Frobenius systems
{\em twist equivalent} if one can be obtained from the other via an
isomorphism and a twist.
Khovanov~\cite{khovanov2} showed that two Frobenius systems related by
twist equivalence give isomorphic link homology groups. It is
important to note however that twisting may ruin nice functoriality
properties with respect to link cobordisms. Actually one can repair things 
again by working with the projective spaces of the homologies, because 
only undesirable scalar factors are caused by twisting.

The following propositions are derived from the work of Lee
\cite{lee}, Shumakovitch \cite{shumakovitch1} and Khovanov
\cite{khovanov2}. For this reason we only sketch the proofs here.

\begin{prop}\label{prop:equivalence}
Let $L$ be a link with $n$ components and let $h,t, \tildeh, \tildet \in \bK$. 
\begin{enumerate}
\item[(i)] If $h^2 + 4t = 0$ then there is an isomorphism
$ \uhtk * L \cong \khk * L$.
\item[(ii)] 
Suppose char$(\bK)\neq 2$. If $h^2 + 4t \neq 0$ and $\frac{\tildeh^2 +
4\tildet}{h^2 +4t} = a^2$ for some non-zero $a\in \bK$ then there is a
twist equivalence $\uhtk * L \cong
\uhtkt * L$.
\end{enumerate}
\end{prop}

\begin{proof}
For (i) let $x$ be the generator of $A_{0,0}$ and $y$ the generator of
$A_{h,t}$. If char$(\bK)\neq 2$ then it can be checked by direct
computation that the map defined by $1\mapsto 1$, $y\mapsto
x+\frac{h}2$ gives an isomorphism of Frobenius systems $A_{h,t} \ra
A_{0,0}$. In characteristic two $h^2+4t=0$ if and only if $h=0$, so
the only non-trivial case is when $t=1$ in which case the map
$1\mapsto 1$, $y\mapsto x+1$ provides an isomorphism.

For (ii) let $x$ be the generator of $A_{h,t}$ and let $y$ be the
generator of $A_{\tildeh,\tildet}$. Let $b=\frac{1}2(\tildeh -ah)$ and
let $A_{h,t}^\prime$ be $A_{h,t}$ twisted by $a^{-1}$. Then by direct
computation one sees that the map $A_{\tildeh,\tildet} \ra
A_{h,t}^\prime$ given by $1\mapsto 1$, $y\mapsto ax+b$ is an
isomorphism of Frobenius systems.
\end{proof}

Note that when $h=0$ and $t=1$ the above result says that
Lee theory over $\FF$ is isomorphic to Khovanov's original theory over
$\FF$, a fact that was proved in \cite{khovanov2}.

\begin{prop}\label{prop:dim}
Let $L$ be a link with $n$ components and
$h,t\in \bK$. If $h^2+4t=\gamma^2$ for some non-zero $\gamma\in \bK$ then
\[
\text{dim}(\uhtk * L) = 2^n.
\] 
All generators lie in even degree and for a knot both generators lie in degree zero.
\end{prop}

\begin{proof}
Change basis to write $A_{h,t} = \bK\{\alpha,\beta\}$ where
\[
\alpha = x-\frac{1}2(h-\gamma)
\]
\[
\beta =  x-\frac{1}2(h+\gamma).
\]
In characteristic two the condition $h^2+4t=\gamma^2 \neq 0$ implies
$h=\gamma=1$, and the basis change is $\alpha = x$ and $\beta = x+1$
which is the change of basis used in \cite{turner}. Courtesy of the condition $h^2+4t=\gamma^2 \neq 0$ this change of
basis diagonalises the multiplication:
\[
\alpha^2 = \gamma\alpha \;\;\;\; \beta^2=-\gamma\beta \;\;\;\; \alpha\beta= \beta\alpha = 0.
\]
The rest of the proof is identical to Lee's proof in \cite{lee} in
which the details of the special case $\bK=\bQ$, $h=0$, $t=1$ and
$\gamma=2$ are provided.
\end{proof}

Khovanov's original link homology was defined integrally and each of the theories discussed so far also has an integral version. Indeed, the
Frobenius system in Example \ref{ex:frobsys} can also be defined over
$\bZ$ resulting in the link homology we denote by $\uhtz * L$.

\begin{prop}\label{prop:int}
Let $L$ be a link with $n$ components and let $h,t\in \bZ$ satisfy
$h^2+4t=\gamma^2$ for non-zero $\gamma\in\bZ$.  
\begin{enumerate}
\item[(i)] There is an isomorphism
\[
\uhtz* L \cong \underbrace{\bZ \oplus \cdots \oplus \bZ}_{2^n} \oplus T^{*}
\]
where $T^{*}$ is all torsion.
\item[(ii)] 
If $h,t<p$ and $\gamma\neq 0$ mod $p$ where $p$ is a prime, then $\uhtz *
L$ has no $p$-torsion.
\end{enumerate}
\end{prop}

\begin{proof}
If $A$ is the Frobenius system giving
$\uhtz * -$ then $A\otimes_\bZ \bQ$ is the Frobenius system giving
$\uhtq * -$. By the construction of link homology this means that each
chain group in the rational theory is the integral chain group
tensored with $\bQ$. Thus the universal coefficient theorem gives
\begin{eqnarray*}
\uhtq i L & \cong &\uhtz i  L \otimes_\bZ \bQ \oplus \Tor^\bZ (\uhtz
     {i+1} L , \bQ )\\
 & = & \uhtz i L \otimes_\bZ \bQ
\end{eqnarray*}
Thus by  Proposition \ref{prop:dim} 
\[
\text{dim}(\uhtz *  L \otimes_\bZ \bQ) = \text{dim}(\uhtq * L) = 2^n
\]
from which part (i) follows.

For part (ii) we will prove by induction on $i$ that $\uhtz i L$ has
no $p$-torsion under the hypotheses given.  Suppose that $\uhtz i L$
has no $p$-torsion for $i\leq N$ and now claim the same holds true for
$i=N+1$. Note that $\uhtz i L$ is non-trivial only for finitely many
values of $i$ so the induction has a base case. By the universal
coefficient theorem we have
\[
\uhtp N L \cong \uhtz N L \otimes_\bZ \bF_p \oplus \Tor^\bZ(\uhtz
     {N+1} L , \bF_p).
\] 
If $N$ is odd, then the left hand side is trivial since it follows from Proposition \ref{prop:dim} that all generators are in even homological degree. Hence
$\Tor^\bZ(\uhtz {N+1} L , \bF_p)= 0$ showing there is no $p$-torsion
in $\uhtz {N+1} L$. If $N$ is even, by Proposition \ref{prop:dim} we
know the number of copies of $\bF_p$ on the left and moreover that the
same number occurs in the first summand on the right, so the Tor group is again trivial and  $\uhtz {N+1}
L$ does not have $p$-torsion.
\end{proof}

For integral Bar-Natan theory one can do slightly better. The change
of basis $\alpha = x$, $\beta=x-1$ in fact diagonalises the theory so
in this case $T^*$ is trivial. For integral Lee theory part (ii) above
shows that the only possible torsion is $2$-torsion.


\section{Rasmussen's $s$-grading}\label{sec:rsg}
As we noted above the theories we are concerned with are not in
general bi-graded but instead possess a filtration. Let $C^*(L)$ be
the complex formed using the Frobenius system $A_{h,t}$ over $\bK$
i.e. whose homology is $\uhtk * L$. As above $\bK$ is one of $\bQ$ or
$\bF_p$ for $p$ a prime and we are assuming $h^2+4t=\gamma^2$ for
$0\neq \gamma \in \bK$.

Define $p\colon C^*(L) \ra \bZ$ as follows. Set $p(1)=1$ and $p(x)=-1$
and for any element $w=w_1\otimes w_2\otimes\cdots\otimes w_m\in
C^*(L)$, where $w_i\in\{1,x\}$, set $p(w)=p(w_1)+ \cdots + p(w_m)$.
An arbitrary $w\in C^*(L)$ is not homogeneous with respect to $p$ but
can be written as $w=w^1+w^2+\cdots w^l$, where $w^j$ is homogeneous
for all $j$. We define 
$$p(w)=\min\left\{p(w^j)\,|\,j=1,\ldots
l\right\}.
$$ 
Now for any $w\in C^{i}(L)$, define
$$
q(w)=p(w)+i+c^+-c^-,
$$ 
where $c^+$ and $c^-$ are the numbers of
positive and negative crossings respectively in $L$. The filtration
grading of an element $w$ is $q(w)$.

As Rasmussen explains in \cite{rasmussen} this determines a grading
$s$ on homology. For $\alpha \in \uhtk * L$ define $$ {\rasleea \alpha
\bK}_{h,t} =\max\left\{q(w)\,|\, w\in C^*(L), [w] = \alpha \right\}.
$$ If there is no confusion we will supress $h$ and $t$ from the
notation writing ${\rasleea \alpha \bK}$ for ${\rasleea \alpha
\bK}_{h,t}$.

For integral theories we define $s(\alpha,\bZ)$ in a similar manner by
restricting the definition to classes $\alpha$ in the torsion-free
part of $\uhtz * L$.

The $s$-grading satisfies some important properties given in the following two propositions.

\begin{prop}\label{prop:sequiv}
Suppose char$(\bK)\neq 2$. If $h^2 + 4t \neq 0$ and $\frac{\tildeh^2 +
4\tildet}{h^2 +4t} = a^2$ for some non-zero $a\in \bK$ then the
twist equivalence of Proposition \ref{prop:equivalence}(ii)
preserves the $s$-grading.
\end{prop}

\begin{proof}
Recall that if $x$ is the generator of $A_{h,t}$ and $y$ is the
generator of $A_{\tildeh,\tildet}$, then the twist equivalence
consists of twisting $A_{h,t}$ by $a^{-1}$ together with an
isomorphism $\psi_*\colon \uhtkt * L \ra \uhtk * L$. This isomorphism
is induced at the level of Frobenius systems by $A_{\tildeh,\tildet}
\ra A_{h,t}^\prime$ defined by $1\mapsto 1$, $y\mapsto ax+b$ where
$b=\frac{1}2(\tildeh -ah)$.

It is clear that the twist preserves $s$ so we only need to consider
the isomorphism above. Let $C^*_{h,t}(L)$ be the complex whose
homology is $\uhtk * L$ and similarly let $C^*_{\tildeh,\tildet}(L)$
be the complex giving $\uhtkt * L$. Let $\psi\colon
C^*_{\tildeh,\tildet}(L) \ra C^*_{h,t}(L)$ be induced by the
isomorphism of Frobenius systems above. We claim that $\psi$ preserves
the filtration degree $q$. We can write $w\in
C^*_{\tildeh,\tildet}(L)$ as
\[
w = \sum\lambda_I\epsilon_I(y)
\]
where each $\epsilon_I(y) = \epsilon_1 \otimes \epsilon_2 \otimes \cdots$
with $\epsilon_j\in \{1,y\}$. By the definition of $\psi$ we have
\[
\psi (\epsilon_I(y)) = a^{r(I)}\epsilon_I(x) + \text{terms of higher filtration}
\]
where $r(I)$ is the number of $y$'s in $\epsilon_I(y)$. From this it
follows that $q(\psi(w)) = q(w)$ since any term $\epsilon_I$ with
$q(\epsilon_I) = q(w)$ also appears in $\psi(w)$.

Next we claim that $\psi_*$ preserves $s$ i.e. for $\alpha\in\uhtkt *
L$
\begin{equation}\label{eq:sleetheta}
{\rasleea \alpha \bK}_{\tildeh,\tildet} = {\rasleea {\psi_*(\alpha)} \bK}_{h,t}.
\end{equation}
Let $w\in C^*_{\tildeh,\tildet}(L) $ such that $[w] = \alpha$ and $q(w) = {\rasleea \alpha \bK}_{\tildeh,\tildet}$. Then $\psi(w)$ represents $\psi_*(\alpha)$ and so
\[
{\rasleea {\psi_*(\alpha)} \bK}_{h,t} \geq q(\psi(w)) = q(w) = {\rasleea \alpha \bK}_{h,t}.
\]
Conversely, let $v\in C^*_{h,t}(L)$ be such that $[v] = \psi_*(\alpha)$ and
$q(v) = {\rasleea {\psi_*(\alpha)} \bK}_{h,t}$.  Then $\psi^{-1}(v)$ represents
$\alpha$ so
\[
{\rasleea \alpha \bK}_{\tildeh,\tildet} \geq q(\psi^{-1}(v)) = q(v) = {\rasleea
	 {\psi_*(\alpha)} \bK}_{h,t}.
\] 
proving (\ref{eq:sleetheta}).
\end{proof}

The next property involves the maps in the universal coefficient theorem. Recall that the universal coefficient theorem provides a short exact sequence
$$
\xymatrix{
0\ar[r] & \uhtz * L\otimes_{\bZ} \bK\ar[r]^-{\phi} & \uhtk * L \ar[r] & \text{Tor}^{\bZ}(\uhtz {*+1} L , \bK) \ar[r] & 0}.
$$

\begin{prop}\label{prop:uct}
If $h^2+4t=\gamma^2$ in $\bZ$ and $\gamma$ is non-zero as an element of $\bK$ then 
\[
\phi\colon \uhtz * L\otimes_{\bZ} \bK \ra \uhtk * L 
\]
is an isomorphism that preserves the $s$-grading.
\end{prop}

\begin{proof}
It is an isomorphism since the Tor group is trivial: over $\bQ$ always
and over $\bF_p$ courtesy of part (ii) of Proposition \ref{prop:int}.

Recall that $\phi$ is induced by the inclusion
\[
\ol{\phi}\colon Z^*(L,\bZ)\otimes \bK \ra C^*(L,\bZ)\otimes \bK = C^*(L,\bK)
\]
which clearly preserves the filtration grading $q$.

To show $\phi$ preserves $s$ we must show that given
$\alpha\in \uhtz * L/\text{Tors}$ we have
\begin{equation}\label{eq:phislee}
\rasleea\alpha \bZ =  
\rasleea{\phi(\alpha\otimes 1)} \bK
\end{equation}

Let $w\in Z^*(L,\bZ)$ be a representative of $\alpha$ such that $q(w)=\rasleea \alpha 
\bZ$. Then $\ol{\phi}(w\otimes 1)$ represents $\phi(\alpha\otimes 1)$
and so
\[
\rasleea{\phi(\alpha\otimes 1)} \bK \geq q(\ol{\phi}(w\otimes 1))=q(w)=\rasleea \alpha 
\bZ. 
\]
Conversely, let $u\in Z^*(L,\bK)$ represent $\phi(\alpha\otimes 1)$
such that $q(u) = \rasleea{\phi(\alpha\otimes 1)} \bK$. 
We may write
$u=\sum v_i \otimes \lambda_i \in Z^*(L,\bZ)\otimes \bK$. When $\bK=\bQ$ let
$\lambda$ be the least common multiple of the denominators of the
$\lambda_i$ and when $\bK=\bF_p$ let $\lambda =1$. Define $v\in Z^*(L,\bZ)$ by
\[
\lambda \sum v_i \otimes \lambda_i = v\otimes 1 \in Z^*(L,\bZ)\otimes \bK.
\]
Note that $q(v) = q(u)$ and moreover that since $\phi$ is an
isomorphism $[v] = \lambda \alpha$. We also have $\rasleea {\lambda\alpha}
\bZ=\rasleea\alpha \bZ$ and so
\[
\rasleea\alpha \bZ = \rasleea {\lambda\alpha} \bZ \geq q(v) = q(u) =
\rasleea{\phi(\alpha\otimes 1)} \bK 
\]
proving (\ref{eq:phislee}) and hence the claim.
\end{proof}


\section{Rasmussen's invariant}\label{sec:ri}
Let $\bK$ be one of $\bQ$ or $\bF_p$ and let $h,t\in\bK$ satisfy $h^2+4t=\gamma^2$ for some $0\neq\gamma\in\bK$. Let $K$ be a knot and define
\[
{\rasleemin \bK}_{h,t} = \min \{ {\rasleea \alpha \bK}_{h,t} \,|\,
\alpha\in \uhtk *K, \alpha\neq 0 \}
\]
and
\[
{\rasleemax \bK}_{h,t} = \max \{ {\rasleea \alpha \bK}_{h,t} \,|\,
\alpha\in \uhtk * K, \alpha\neq 0 \}.
\]

Rasmussen's $s$-invariant for the theory $\uhtk * -$ is defined as follows. The original definition in \cite{rasmussen} is for the case $\bK=\bQ$. 

\begin{defn}
\[
{\rasleeknot\bK}_{h,t}=\dfrac{{\rasleemin\bK}_{h,t}+{\rasleemax\bK}_{h,t}}{2}
\]
\end{defn}

For integral theories we may make an analogous definition by using
$s(\alpha,\bZ)$ which we recall restricts its definition to the the
torsion-free part of $\uhtz * K$.

Here is our main result. Let $\bK$ and $\widetilde \bK$ be $\bQ$ or
$\bF_p$ ($\bK$ and $\widetilde\bK$ possibly different).

\begin{thm}\label{thm:main}
Let $K$ be a knot. Let $h,t,\tilde h,\tilde t\in \bZ$ be such that
$h^2+4t=\gamma^2\ne 0$ and $\tildeh^2+\tildet=\tilde \gamma^2\ne 0$
with $\gamma\ne 0\in \bK$ and $\tilde\gamma\ne 0\in
\widetilde\bK$. Then
\[
{\rasleeknot\bK}_{h,t} = {\rasleeknot{\widetilde\bK}}_{\tilde h,\tilde t}
\]
holds.
\end{thm}

\begin{proof}
Using Proposition \ref{prop:equivalence}(ii) and Proposition 
\ref{prop:sequiv} we have 
\begin{equation}\label{eq:part1}
{\rasleeknot\bQ}_{h,t} = {\rasleeknot\bQ}_{\tilde h,\tilde t}.
\end{equation}
Letting $\bK^\prime$ be any of $\bQ$ or $\bF_p$, Proposition \ref{prop:uct}
implies 
\begin{equation}\label{eq:part2}
{\rasleeknot\bKprime}_{h,t} = {\rasleeknot\bZ}_{h,t}.
\end{equation}
From (\ref{eq:part1}) and (\ref{eq:part2}) it follows that 
\[
{\rasleeknot\bK}_{h,t} = {\rasleeknot\bZ}_{h,t} = {\rasleeknot\bQ}_{h,t} =
{\rasleeknot\bQ}_{\tilde h,\tilde t} =  {\rasleeknot\bZ}_{\tilde h,\tilde t} =
  {\rasleeknot{\widetilde\bK}}_{\tilde h,\tilde t}.
\]
\end{proof}

In particular ${\rasleeknot\FF}_{1,0} = {\rasleeknot\bQ}_{0,1}$ showing that 
the $s$-invariant from Bar-Natan's characteristic two theory is equal to 
Rasmussen's original $s$-invariant defined using Lee theory over $\bQ$.

\vspace*{1cm}

\noindent {\bf Acknowledgements} The first author was supported by the 
Funda\c {c}\~{a}o para a Ci\^{e}ncia e a Tecnologia through the
programme ``Programa Operacional Ci\^{e}ncia, Tecnologia, Inova\c
{c}\~{a}o'' (POCTI), cofinanced by the European Community fund FEDER.

The second author was supported by the
European Commission through a Marie Curie fellowship and thanks the
Institut de Recherche Math\'ematiques Avanc\'ee in Strasbourg for
their hospitality.



\begin{thebibliography}{15}
\bibitem{barnatan2}
D. Bar-Natan, Khovanov's homology for tangles and cobordisms, 
{\em Geometry and Topology}, 9 (2005), paper no. 33, 1443-1499.
\bibitem{khovanov1}
M. Khovanov, A categorification of the Jones polynomial, {\em Duke
Math J.}, 101 (2000), 359-426.
\bibitem{khovanov2} 
M. Khovanov, Link homology and Frobenius extensions, math.QA/0411447, 
2004.
\bibitem{lee}
E. Lee, An endomorphism of the Khovanov invariant,
math.GT/0210213, 2002.
\bibitem{rasmussen}
J. Rasmussen, Khovanov homology and the slice genus, math.GT/0402131,
2004.
\bibitem{rasmussensurvey}
J. Rasmussen, Knot polynomials and knot homologies, math.GT/0504045, 2005.
\bibitem{shumakovitch1}
A. Shumakovitch, Torsion of the Khovanov homology, math.GT/0405474, 2004.
\bibitem{turner}
P. Turner, Calculating Bar-Natan's characteristic two Khovanov
homology, math.GT/0411225, 2004.
\end{thebibliography}
\end{document}